\documentclass[12pt,reqno]{amsart}
\usepackage{amsmath,amsthm,amssymb,amscd,url,enumerate}
\usepackage{graphicx}
\usepackage{afterpage}
\usepackage{tikz}
\usepackage[all]{xy}
\usepackage[colorlinks=true]{hyperref}
\usepackage[width=5.3in,height=8.5in, centering]{geometry}

\newcommand{\TITLE}{An arborist's guide to the rationals}
\newcommand{\TITLERUNNING}{An arborist's guide to the rationals}
\newcommand{\DATE}{\today}
\newcommand{\VERSION}{1}

\theoremstyle{plain} 

\newtheorem*{theorem*}{Theorem}

\theoremstyle{definition}
\newtheorem*{definition}{Definition}

\theoremstyle{remark}

\numberwithin{theorem}{section}

  {\end{list}}

  {\end{list}}

\newcommand{\tightoverset}[2]{%
  \mathop{#2}\limits^{\vbox to -.5ex{\kern-1.05ex\hbox{$#1$}\vss}}}

\newcommand{\PP}{\mathbb{P}}

\newcommand{\QQ}{\mathbb{Q}}

\newcommand{\ZZ}{\mathbb{Z}}

\newcommand\Conway{MR1478672}
\newcommand\CalkinWilf{MR1763062}

\newcommand\Bates{MR2673006}
\newcommand\Glasby{MR2854004}
\newcommand\Guthery{MR2895290}
\newcommand\LagariasTresser{MR2361372}

\newcommand\Reznick{MR1084197}
\newcommand\Lansing{Lansing}
\newcommand\JoaoOne{MR2672845}
\newcommand\JoaoTwo{MR2789032}
\newcommand\BirdTree{MR2539335}

\title[\TITLERUNNING]{\vspace*{-1.3cm} \TITLE}

\author{Katherine E. Stange}

\date{\DATE, Draft \#\VERSION}
\address{%
Department of Mathematics, University of Colorado,
Campux Box 395, Boulder, Colorado 80309-0395}
\email{kstange@math.colorado.edu}

\thanks{The author's research has been supported by NSA Grant Number H98230-14-1-0106
}

\begin{document}


\begin{abstract}
        There are two well-known ways to enumerate the positive rational numbers in an infinite binary tree:  the Farey/Stern-Brocot tree and the Calkin-Wilf tree.  In this brief note, we describe these two trees as `transpose shadows' of a tree of matrices (a result due to Backhouse and Ferreira) via a new proof using yet another famous tree of rationals:  the topograph of Conway and Fung.
\end{abstract}

\maketitle

\section{Four Trees}

In 2000, Calkin and Wilf studied an explicit enumeration of the positive rationals which naturally arranges itself into an infinite tree \cite{\CalkinWilf}, the first few levels of which are shown here:

\begin{center}
\includegraphics[width=4in]{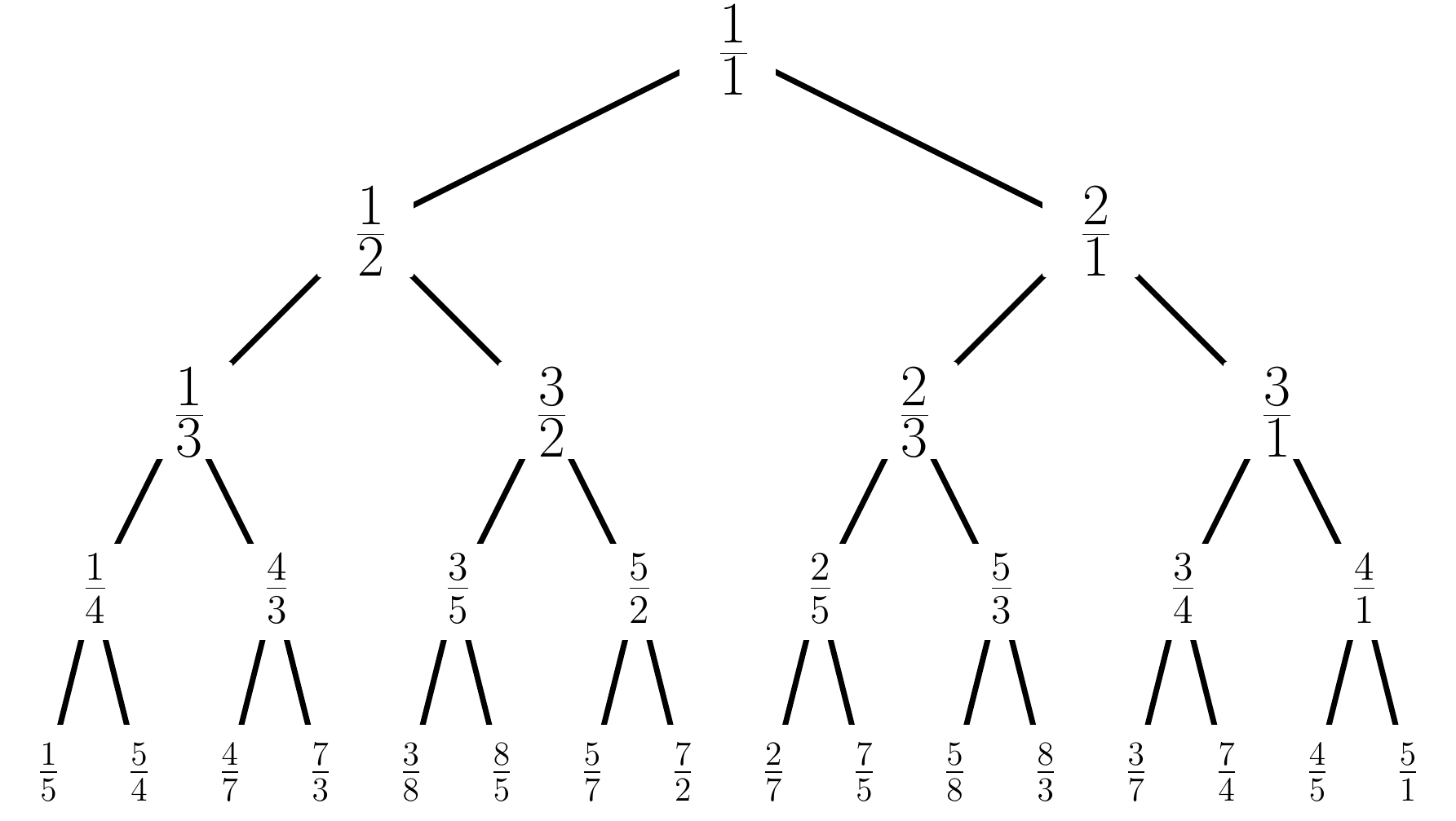}
\end{center}

The generation rule is that a parent $\frac{a}{b}$ has the following {left} and {right} children:

\begin{center}
\includegraphics[width=1in]{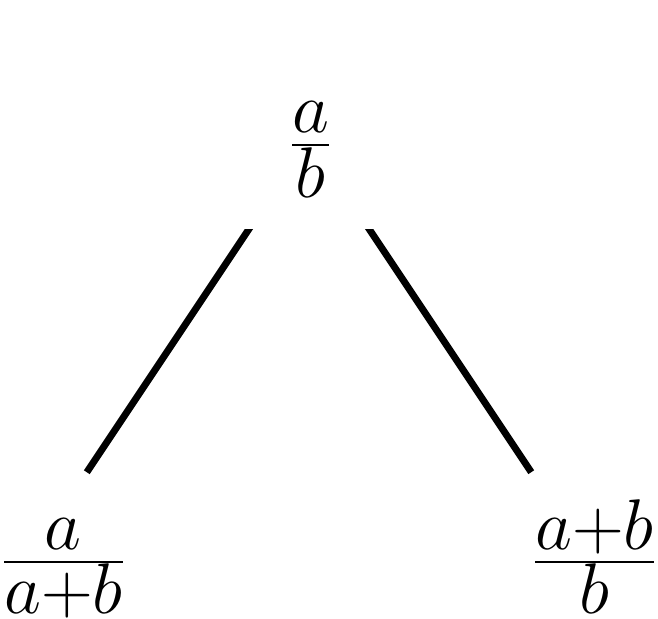}
\end{center}

Every positive rational number appears in this tree exactly once.  Calkin and Wilf consider the integer sequence $b(n)$ which enumerates the representations of $n$ as a sum of powers of $2$, where each power is allowed to appear at most twice.  The function $b(n)/b(n~+~1)$ reads off the entries in the tree left to right, top row downwards:
\[
        \frac{1}{1}, \frac{1}{2}, \frac{2}{1}, \frac{1}{3}, \frac{3}{2}, \frac{2}{3}, \frac{3}{1}, \frac{1}{4}, \frac{4}{3}, \frac{3}{5}, \ldots
\]

This tree is reminiscent of the more famous Farey tree, also known as the Stern-Brocot tree, which begins as follows:

\begin{center}
\includegraphics[width=4in]{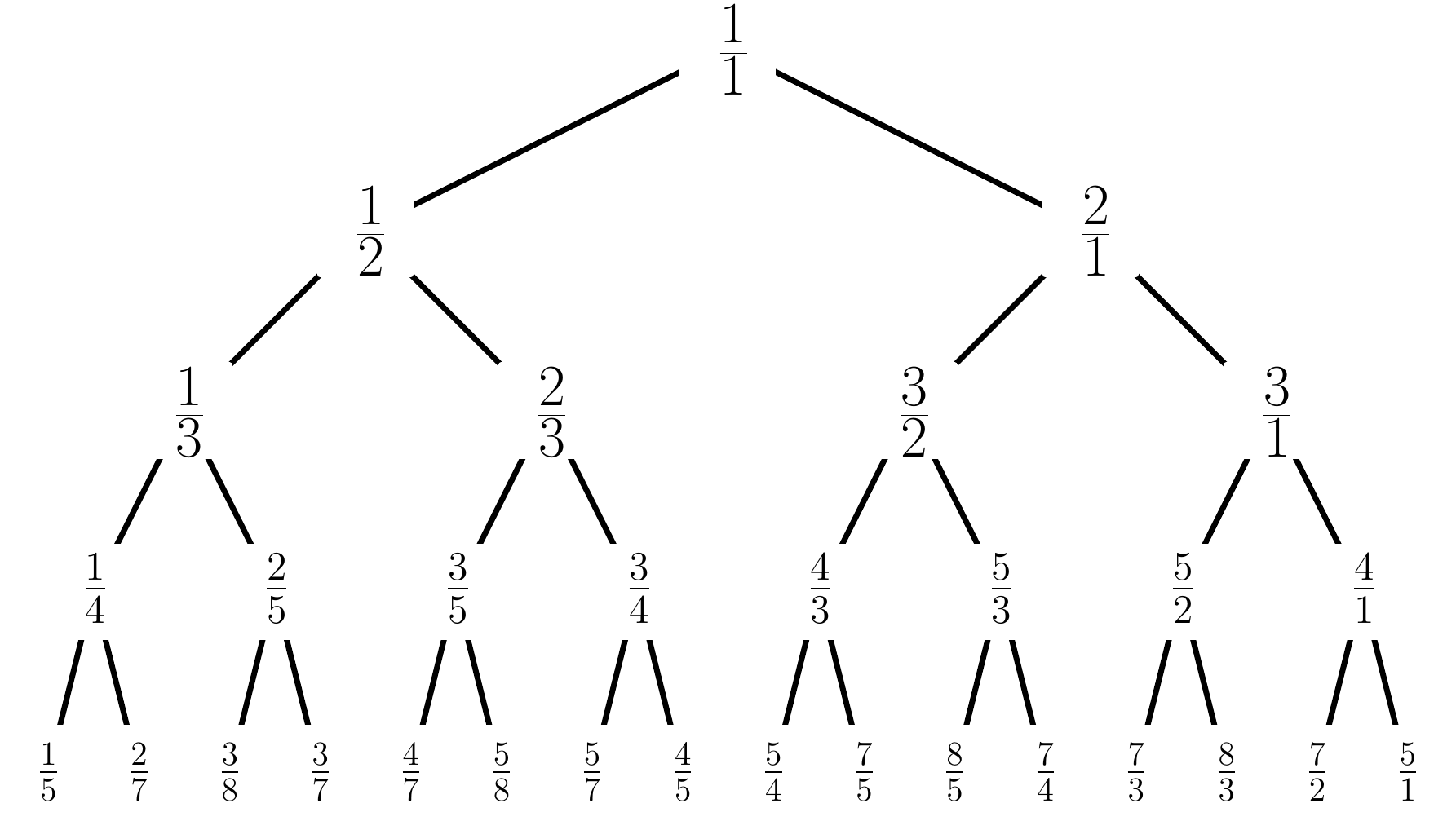}
\end{center}

The latter name is in honour of two independent descriptions of related ideas in the mid 1800's by Stern \cite{Stern} and Brocot \cite{Brocot}.  Brocot was a french clockmaker who created an array of fractions for the purpose of designing clockwork gears\footnote{Brocot wrote a book and a paper by the same title, \emph{Calcul des rouages par approximation, nouvelle m\'ethode}; it is the book which contains the array.}.  Stern studied an array of integers, which can be used to generate \emph{both} the Stern-Brocot tree and the Calkin-Wilf tree (it has been quite reasonably suggested the latter tree be called the \emph{Eisenstein-Stern tree}, but this name is not prevalent \cite{\JoaoTwo}).   The name `Farey tree' comes from its relationship to Farey sequences\footnote{Given a bound $D$, the associated Farey sequence is the sequence of rationals in $[0,1]$ which, in lowest form, have denominator less than or equal to $D$.  The name has gradually become associated to a wider variety of structures generated by means of the mediant operation we will describe in a moment.} (which were themselves likely invented by Charles Haros \cite{\Guthery}).  For more on the muddy historical waters, and the trees themselves, see \cite{\JoaoTwo,\LagariasTresser,\Lansing}.

The \emph{mediant} of rationals $\frac{a}{b}$ and $\frac{c}{d}$ (in lowest form) is $\frac{a+c}{b+d}$.  The root of the tree is $\frac{1}{1}$, which forms the first row in the tree.  Bracket this row by $\frac{0}{1}$ and $\frac{1}{0}$, and then take the list of mediants:
\vspace*{-1em}

\[
        \xymatrix @R=0.15pc @C=0.40pc {
                \frac{0}{1} \ar[dr] & & \frac{1}{1} \ar[dl] \ar[dr] & & \frac{1}{0} \ar[dl] \\
                             & \frac{1}{2} & & \frac{2}{1} & \\
        }
\]
\vspace*{-1em}

These mediants form the second row of the tree.  In general, beginning with the full list of fractions appearing in rows $1$ through $n$, listed in order of size, ones brackets as above, resulting in what has been called a \emph{Brocot sequence} or a \emph{Farey-like sequence}:
\[
        \frac{0}{1}, \frac{1}{n}, \ldots, \frac{n}{1}, \frac{1}{0}.  
\]
The mediants of the Brocot sequence form the $(n+1)$-st row.  Continue ad infinitum and the tree will, just as the Calkin-Wilf tree does, contain exactly one instance of each positive rational number.  

It should not be surprising that these two trees share a common genesis.  Stern's array gives rise to the \emph{Stern sequence} $s(n)$.  The $b(n)$ of the Calkin-Wilf tree is exactly $s(n+1)$ \cite{\Reznick}, while the fractions of the Farey tree are of the form $s(n)/s(2^r-n)$ (see \cite{\Lansing} for history).  Some algebraic connections between these two trees are described in \cite{\Bates, \Glasby}.  

It is the purpose of our arboreal tour to explore a connection between these two famous rational-enumerating trees via a single tree of \emph{matrices}.
\tiny
\[
        \xymatrix @R=1.25pc @C=-1.40pc {
                & & & & & & & &  {\begin{pmatrix} 1 & 0 \\ 0 & 1 \end{pmatrix}} \ar@{-}[lllld] \ar@{-}[rrrrd] & & & & & & & & \\
                & & & & {\begin{pmatrix} 1 & 0 \\ 1 & 1 \end{pmatrix}}  \ar@{-}[lld] \ar@{-}[rrd] & & & & & & & & {\begin{pmatrix} 1 & 1 \\ 0 & 1 \end{pmatrix}}   \ar@{-}[lld] \ar@{-}[rrd] & & & & & & & & & \\
                & &
                {\begin{pmatrix} 1 & 0 \\ 2 & 1 \end{pmatrix}} 
                \ar@{-}[ld] \ar@{-}[rd] 
                & & & &
                {\begin{pmatrix} 2 & 1 \\ 1 & 1 \end{pmatrix}}
                \ar@{-}[ld] \ar@{-}[rd] 
                                   & & & &
                {\begin{pmatrix} 1 & 1 \\ 1 & 2 \end{pmatrix}} 
                \ar@{-}[ld] \ar@{-}[rd] 
                                   & & & &
                {\begin{pmatrix} 1 & 2 \\ 0 & 1 \end{pmatrix}}
                \ar@{-}[ld] \ar@{-}[rd] 
                & & & & &\\
                & {\begin{pmatrix} 1 & 0 \\ 3 & 1 \end{pmatrix}} & &  
                {\begin{pmatrix} 3 & 1 \\ 2 & 1 \end{pmatrix}} & &  
                {\begin{pmatrix} 2 & 1 \\ 3 & 2 \end{pmatrix}} & &   
                {\begin{pmatrix} 3 & 2 \\ 1 & 1 \end{pmatrix}} & &    
                {\begin{pmatrix} 1 & 1 \\ 2 & 3 \end{pmatrix}} & &  
                {\begin{pmatrix} 2 & 3 \\ 1 & 2 \end{pmatrix}} & &  
                {\begin{pmatrix} 1 & 2 \\ 1 & 3 \end{pmatrix}} & &   
                {\begin{pmatrix} 1 & 3 \\ 0 & 1 \end{pmatrix}} & & \\
        }
\]
\normalsize

To obtain this tree, place the $2 \times 2$ identity matrix at the root, and apply the following generation rule:
\small
\[
        \xymatrix @!=1.7pc {
                 & & M \ar@{-}[ld] \ar@{-}[rd] & &    \\
                &  {\begin{pmatrix} 1 & 0 \\ 1 & 1 \end{pmatrix}}M   &  & {\begin{pmatrix} 1 & 1 \\ 0 & 1 \end{pmatrix}}M   &   \\
        }
\]
\normalsize

The matrix tree is a visualization of the folk theorem that the monoid $\operatorname{SL}_2(\ZZ^{\ge 0})$ is freely generated by the two elements
\[
        \begin{pmatrix} 1 & 0 \\ 1 & 1 \end{pmatrix}, \quad \mbox{and}\quad
        \begin{pmatrix} 1 & 1 \\ 0 & 1 \end{pmatrix};
\]
every element of $\operatorname{SL}_2(\ZZ^{\ge 0})$ appears in the tree exactly once.  See Section \ref{sec:mob}.

The relationship between these three trees is originally due to Backhouse and Ferreira, and deserves to be better known.

\begin{theorem*}[Backhouse, Ferreira \cite{\JoaoOne, \JoaoTwo}]
To recover the Calkin-Wilf tree, one replaces, in the matrix tree above,
\[
        \begin{pmatrix} a & b \\ c & d \end{pmatrix} \quad \mbox{with} \quad
        \frac{a+b}{c+d}.
\]
To recover the Farey tree, one replaces
\[
        \begin{pmatrix} a & b \\ c & d \end{pmatrix}
        \quad \mbox{with} \quad \frac{d+b}{c+a}.
\]
\end{theorem*}

The relation to the Calkin-Wilf tree is immediate by comparing the generation rules of the two trees, and this relationship was exploited in \cite{Kucharczyk, Nathanson}.

The relationship to the Farey tree is also not too difficult to verify directly, but it is our purpose to provide a new proof which arises by turning to yet another beautiful tree that enumerates the rationals: the topograph of Conway and Fung \cite{\Conway}.

\begin{center}
\includegraphics[width=4in]{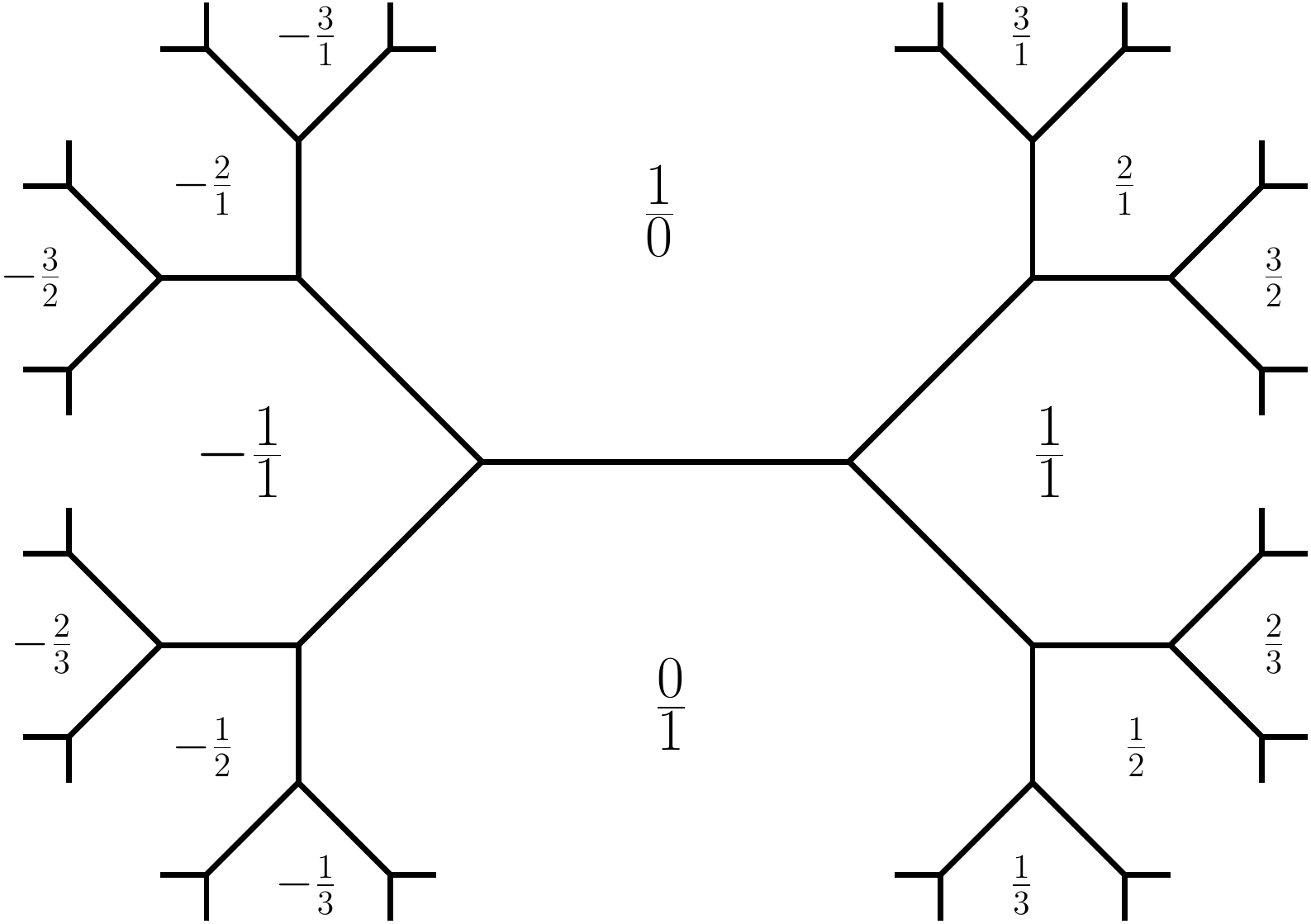}
\end{center}

This time, it is the regions that are labelled by the rational numbers.  Surrounding each vertex are three fractions which are, in some order (and with appropriate use of signs), a pair of fractions together with their mediant.

{\bf Acknowledgements.}  The author would like to thank Jo\~ao Ferreira for gently bringing to her attention that the result was previously known.  She would also like to thank Jo\~ao Ferreira, Jeffrey Lagarias and Jennifer Lansing for help in correcting the historical background contained in this article.  Thank you also to Jesse Levine and Amy Feaver for comments on a draft.

\section{The topograph}

Write $\infty = \frac{1}{0}$, and $\QQ^\infty = \QQ \cup \{ \infty \}$.

\begin{definition}
        Two points $\frac{a}{b}, \frac{c}{d} \in \QQ^\infty$ (given in lowest terms) are called \emph{$\ZZ$-distinct} if $ad-bc = \pm 1$.
\end{definition}

The definition is symmetric and doesn't depend on the convention for minus signs in one's definition of `lowest terms.'

\begin{definition}
The \emph{topograph} is the graph whose set of vertices is all triples of pairwise $\ZZ$-distinct points, with the stipulation that two such triples are connected by an edge whenever they have a pair of elements in common.  
\end{definition}

We can identify an edge with the pair of $\ZZ$-distinct elements shared by its vertices\footnote{Conway and Fung consider $\QQ^\infty$ as $\PP^1(\QQ)$, so that points become primitive vectors of $\ZZ^2$; each edge corresponds to a basis of $\ZZ^2$, and vertices give triples called \emph{superbases}.}.  Note that every pair appears exactly once in the topograph since the pair ${\frac{a}{b},\frac{c}{d}}$ is part of only two triples:
\[
        \left\{ \frac{a}{b}, \frac{c}{d}, \frac{a+c}{b+d} \right\}, \quad \mbox{and}
        \quad
        \left\{ \frac{a}{b}, \frac{c}{d}, \frac{a-c}{b-d} \right\}.
\]

The graph can be made planar in such a way that the boundary of each region (an infinite tree of valence 2, i.e. a line), consists of all pairs and triples containing a fixed element.  In this way, each region can be labelled with a point of $\QQ^\infty$ \cite{\Conway}.

\section{M\"obius transformations}
\label{sec:mob}

The automorphisms of $\QQ^\infty$ are the M\"obius transformations,
\[
        z \mapsto \frac{az+b}{cz+d}, \quad a,b,c,d \in \QQ, \quad ad-bc \neq 0,
\]
forming a group under composition.  This is isomorphic to the matrix group
\[
        \operatorname{PGL}_2(\QQ) = \left\{ \begin{pmatrix} a & b \\ c & d \end{pmatrix} : a,b,c,d \in \QQ, ad-bc \in \QQ^* \right\} / \left\{ k I_{2\times 2}: k \in \QQ^* \right\}, 
\]
by the map
\[
        \left(z \mapsto \frac{az+b}{cz+d}\right) \mapsto \begin{pmatrix} a & b \\ c & d \end{pmatrix}.
\]
The subset of matrices having representatives with non-negative integer entries and determinant $1$ is closed under multiplication but not inverses, forming the monoid 
\[
        \operatorname{SL}_2(\ZZ^{\ge 0}) = \left\{ \begin{pmatrix} a & b \\ c & d \end{pmatrix}: a,b,c,d \in \ZZ^{\ge 0},ad-bc=1 \right\}.
\]
The monoid $\operatorname{SL}_2(\ZZ^{\ge 0})$ is closed under transposition, 
\[
        \gamma \mapsto \gamma^T, \quad
        \frac{az+b}{cz+d} \mapsto
        \frac{az + c}{bz + d}.
\]
With this notation, the Theorem can now be phrased as follows: replacing the transformation $\gamma$ with $\gamma(1)$ gives the Calkin-Wilf tree, while replacing it with $1/\gamma^T(1)$ gives the Farey tree\footnote{The reciprocal is immaterial, since it would disappear if we wrote the Farey tree right-to-left instead of left-to-right, i.e. reflected in its vertical midline.}.

\section{The Proof}

The proof proceeds by labelling the topograph two ways:  first, to create the Farey tree, and second, to create the matrix tree.  Comparing the two labellings generates the rule given in the Theorem.

\begin{definition}
        A \emph{flow} of the topograph (or a portion of it) is an assignment of direction to every edge in such a way that in-degree is exactly one at each vertex.
\end{definition}

Once one edge is assigned a direction, the portion of the topograph that is forward of that edge (according to the edge direction), has all its directions determined uniquely by the condition of flow, and forms a rooted binary tree directed away from the root.  Choosing the edge $\{ 0, \infty \}$, and directing it toward the vertex $\{ 0, \infty, 1\}$, we obtain the following.

\begin{center}
\includegraphics[width=4in]{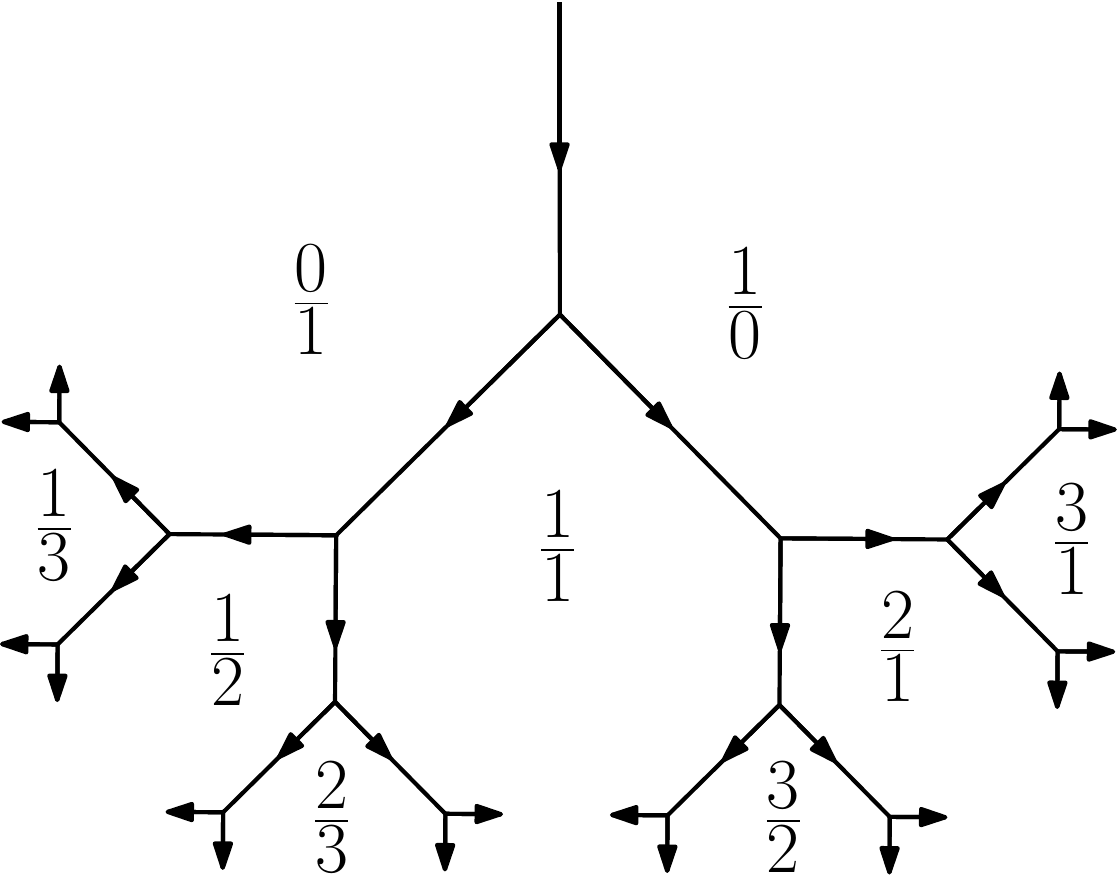}
\end{center}

Each vertex has one incoming edge; with respect to this direction, there's a left, right and forward region.  If we label a vertex with the region bounded by the two outgoing edges (i.e. moving the region labels up to the `peaks' of their respective regions), we obtain the Farey tree \cite{\Conway}.  In particular, all regions are labelled with positive rationals.

By contrast, to such a vertex we may also associate a M\"obius transformation $\gamma(z)$ by specifying 
its values 
at $\frac{1}{0}$, $\frac{0}{1}$ and $\frac{1}{1}$ are exactly the labels of the regions to the left, right and below the vertex.  
In other words, if these labels are, respectively, $\frac{a}{c}$, $\frac{b}{d}$, and $\frac{a+c}{b+d}$, then we obtain the transformation \tiny $\begin{pmatrix} a & b \\ c & d \end{pmatrix}$\normalsize.

\begin{center}
\begin{minipage}{.45\textwidth}
\end{minipage}%
\begin{minipage}{.25\textwidth}
  \includegraphics[width=1in]{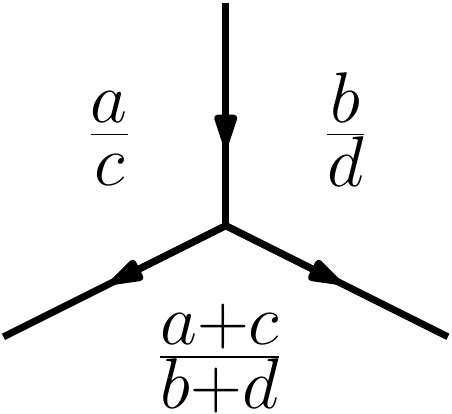}
\end{minipage}%
\begin{minipage}{.30\textwidth}
$\mapsto\quad\quad \begin{pmatrix} a & b \\ c & d \end{pmatrix}$
\end{minipage}
\end{center}

This gives a tree of matrices.  In fact, this is \emph{not} the matrix tree in the introduction, but by applying the map
\[
        \gamma(z) \mapsto (1/\gamma(z))^T, \quad
        \begin{pmatrix} a & b \\ c & d \end{pmatrix} \mapsto
        \begin{pmatrix} c & a \\ d & b \end{pmatrix},
\]
at each vertex, we obtain the matrix tree of the introduction.  The verification of this is straightforward: the roots agree and the $\ZZ$-distinctness condition at each vertex translates into the matrix tree generation rule.  This is sufficient to complete the proof.

\section*{Afterthoughts}

There is another rational-enumerating tree, the Bird Tree \cite{\BirdTree}, which is a levelwise permutation of the Farey tree, but sadly we won't explore it in this note.  It is also interesting to observe that if we define a flow of the topograph which directs the boundary of the region $\frac{1}{0}$ from negative to positive regions, we obtain the extended Farey tree of \cite{\LagariasTresser}.

\bibliography{calkin-wilf}{}
\bibliographystyle{plain}

\end{document}